\documentclass[12pt]{amsart}
\usepackage{amssymb}
\usepackage{epsfig}
\usepackage{graphicx}
\numberwithin{equation}{section}

\input xy
\xyoption{all}

\advance\headheight 2pt
\calclayout
\allowdisplaybreaks[3]

\theoremstyle{plain}
\newtheorem{prop}{Proposition}

\newtheorem{theo}[prop]{Theorem}
\newtheorem{coro}[prop]{Corollary}

\newtheorem{lemm}[prop]{Lemma}

\theoremstyle{definition}
\newtheorem{defi}[prop]{Definition}

\newtheorem{rema}[prop]{Remark}

\newtheorem{exam}[prop]{Example}

\newcommand{\bC}{\mathbb C}

\newcommand{\bP}{\mathbb P}
\newcommand{\bQ}{\mathbb Q}
\newcommand{\bZ}{\mathbb Z}

\newcommand{\cD}{\mathcal D}
\newcommand{\cI}{\mathcal I}
\newcommand{\cK}{\mathcal K}
\newcommand{\cL}{\mathcal L}

\newcommand{\cN}{\mathcal N}
\newcommand{\cO}{\mathcal O}
\newcommand{\cS}{\mathcal S}
\newcommand{\cT}{\mathcal T}

\newcommand{\cX}{\mathcal X}

\newcommand{\ra}{\rightarrow}

\newcommand{\Spec}{{\rm Spec}}
\newcommand{\Gr}{{\rm Gr}}
\newcommand{\NS}{{\rm NS}}
\newcommand{\Hom}{{\rm Hom}}
\newcommand{\iHom}{{\mathcal H}\!{\mathit om}}
\newcommand{\Hilb}{{\mathcal H}\!{\mathit ilb}}
\newcommand{\Ext}{{\rm Ext}}
\newcommand{\iExt}{{\mathcal E}\!{\mathit xt}}

\makeatother
\makeatletter

\author{Brendan Hassett}
\address{Department of Mathematics \\
Rice University, MS 136 \\
Houston, TX 77251-1892}
\email{hassett@rice.edu}

\author{Yuri Tschinkel}
\address{
                Mathematisches Institut \\
                Bunsenstr. 3-5  \\
        37073 G\"ottingen \\
                Germany }
\email{
yuri@uni-math.gwdg.de}

\title[Potential density of rational points]{Potential
density of rational points for K3 surfaces
over function fields}

\begin{document}
\date{\today}

\begin{abstract}
We give examples of non-isotrivial
K3 surfaces over complex function fields with Zariski-dense rational
points and N\'eron-Severi rank one.  
\end{abstract}

\maketitle

\section{Introduction}
\label{sect:introduction}

Let $F$ be a number field or a function field of a curve over the complex numbers.  
For each variety $X$ defined over $F$, we can ask whether the rational points
$X(F)$ are dense in $X$.  We say that rational points of $X$ are {\em potentially dense}
if there exists a finite extension $E/F$ with $X(E)$ dense in $X$.  Potential
density of rational points is expected to be a geometric property, depending (for
smooth projective $X$) on the positivity of the canonical class $K_X$.  For example,
if $X$ is a curve or a surface with $K_X$ negative, rational points are potentially dense;
indeed, such varieties are known to be rational over 
some finite extension of $F$.
On the other hand, potential density fails for curves and for certain classes
of higher-dimensional varieties with $K_X$ positive, by a theorem of Faltings.  
We refer the reader to \cite{Ca} for a general discussion 
of the connections between classification theory and potential 
density questions.

The intermediate case where $K_X=0$ remains mysterious.  Already the case of K3 surfaces
is open, in general.  Potential density has been proved for abelian varieties, Enriques
surfaces, and special classes of K3 surfaces, e.g., those with an elliptic fibration
or an infinite group of automorphisms \cite{HarTs} \cite{bt} \cite{bt1} \cite{h}.  However, all the 
existing examples of K3 surfaces with dense rational points have geometric
N\'eron-Severi rank
at least two.  

Our main result is:
\begin{theo}
\label{theo:one}
Let $B$ be a complex curve and $F=\bC(B)$ its function field. 
There exist non-isotrivial K3 surfaces over $F$, with 
geometric N\'eron-Severi rank one,
and a Zariski dense set of $F$-rational points.
\end{theo}
This partially answers a question of Campana \cite[Sect. 11]{Ca}.

Over a number field, no examples of K3 surfaces with geometric
N\'eron-Severi rank one and a dense set of 
rational points are known.  Our method uses the fact that the ground field is uncountable (see
Remark~\ref{rema:uncountable});  it does not yield
examples over $\bar{\bQ}(B)$.  

\

We summarize the contents of this paper.  Sections~\ref{sect:deform}
through \ref{sect:fe} contain the proof of Theorem~\ref{theo:one}.
The reader may want to consult the heuristic discussion in
Section~\ref{sect:proj} before tackling the main argument.
Our general strategy is to deform reducible curves--unions of 
sections with smooth
rational curves in the fibers--to new sections of larger degree
(see Proposition~\ref{prop:oneps}).  
Sections~\ref{sect:deform} and \ref{sect:k3} review general results
on deformations of curves and moduli spaces of K3 surfaces.
Our examples are pencils of K3 surfaces of low degree;
Section~\ref{sect:k3pencil} describes these pencils
and Proposition~\ref{prop:makepencil} indicates the situations where
our construction applies.  It is an open problem to find examples
of large degree with the properties described in Theorem~\ref{theo:one}. 
Indeed, recent results of Gritsenko, Hulek, and Sankaran
imply that the moduli space of K3 surfaces of large degree is
of general type;  in particular, there are no pencils through
the generic point.  Section~\ref{sect:limit} partially explains 
why pencils play a special r\^ole in our analysis.  Section~\ref{sect:CY}
extends our results to higher-dimensional Calabi-Yau hypersurfaces.

\

\noindent {\bf Acknowledgments:}  The first author was supported by the Sloan Foundation and 
NSF grant 0134259.  We benefited from conversations with
Izzet Coskun, Aise Johan de Jong, Joe Harris, and Jason Starr.  
We are grateful to Coskun and Starr for allowing us to include
their results in Section \ref{sect:CY}.

\section{Deformation theory}
\label{sect:deform}

In this section we recall basic facts from deformation theory used in subsequent proofs.
Throughout, we work over an algebraically closed field $k$.

Let $X$ be a smooth projective variety over $k$, $C\subset X$ a closed
subscheme.  The tangent space
to the Hilbert scheme $\Hilb(X)$ at $[C]$ is 
$$
\cT_{\Hilb(X)}|_{[C]}=\Hom_X(\cI_C,\cO_C)=\Hom_C(\cI_C/\cI_C^2,\cO_C).
$$
We have an exact sequence
$$
\cI_C/\cI^2_C\ra \Omega^1_{X}\otimes_{\cO_X} \cO_C \ra \Omega^1_C\ra 0;
$$
when $C$ is reduced, the left arrow is injective at generic points of $C$ 
\cite{Mat} pp. 188, \cite{Hart} II.8 Ex.1.  

Assume now that $C$ is a reduced local-complete-intersection
subscheme of $X$, so that $\cI_C/\cI_C^2$ 
is a locally free $\cO_C$-module.  Since $C$ is reduced,
an $\cO_C$-linear homomorphism of locally-free modules is injective provided
it is injective at each generic point, and we obtain
$$
0 \ra \cI_C/\cI^2_C\ra \Omega^1_{X}\otimes_{\cO_X} \cO_C \ra \Omega^1_C\ra 0.
$$
Dualizing our exact sequence yields
\begin{eqnarray*}
& & 0 \ra  \iHom_C(\Omega^1_C,\cO_C) \ra \iHom_C(\Omega^1_X,\cO_C) \ra 
\iHom_C(\cI_C/\cI_C^2,\cO_C)  \quad \quad (\dagger) \\
& & \quad \ra \iExt^1_C(\Omega^1_C,\cO_C)\ra 0 
\end{eqnarray*}
on the level of coherent sheaves and
\begin{eqnarray*}
& & 0 \ra \Hom_C(\Omega^1_C,\cO_C) \ra \Hom_C(\Omega^1_X,\cO_C) \ra \Hom_C(\cI_C/\cI^2_C,\cO_C)  \\
& & \quad \ra \Ext^1_C(\Omega^1_C,\cO_C)  
\end{eqnarray*}
on cohomology.  We write
$$\cT_C=\iHom_C(\Omega^1_C,\cO_C) \
\cT_X|_C=\iHom_C(\Omega^1_X,\cO_C) \
\cN_{C/X}=\iHom_C(\cI_C/\cI_C^2,\cO_C);
$$
the last two of these are locally free.  We can rewrite $(\dagger)$ as
\begin{equation}
\label{ES2}
0 \ra \cT_C \ra \cT_X|_C \ra \cN_{C/X} \ra \iExt^1_C(\Omega^1_C,\cO_C) \ra 0.
\end{equation}
Moreover, we have
$$\cT_{\Hilb(X)}|_{[C]}=\Gamma(C,\cN_{C/X}).$$

Assume that $C$ is a nodal curve and fix $p\in C$ a node.  We can choose
local coordinates $x_1,\ldots,x_n$ for $X$ centered about $p$ so that
$$C=\{x_1x_2=x_3=\ldots=x_n=0 \}.$$
We give local presentations for the terms in exact sequence (\ref{ES2}).
The tangent bundle of $X$ restricted to $C$ is freely generated by
$$\partial/\partial x_1, \partial/\partial x_2,\ldots, \partial/\partial x_n.$$
The subsheaf $\cT_C$ is generated by
$ x_1\partial/\partial x_1-x_2\partial/\partial x_2 $
which induces the relation
$$ x_1\partial/\partial x_1=x_2\partial/\partial x_2 $$
in the quotient sheaf $\cT_X|_C/T_C$.  The sheaf $\iExt^1_C(\Omega^1_C,\cO_C)$
has fiber at $p$ equal to the versal deformation space of the isolated singularity
$(C,p)\subset X$, which in this case is one-dimensional.  The normal bundle $\cN_{C/X}$
is an extension of this $\Ext$-group by the quotient $\cT_X|_C/T_C$, with free
generators
\begin{equation} \label{GEN}
x_2^{-1}\partial/\partial x_1=x_1^{-1}\partial/\partial x_2,
\partial/\partial x_3,\ldots, \partial/\partial x_n. \end{equation}

Let $\nu:C^{\nu}\ra C$ denote the normalization of $C$ at the point $p$,
with conductor $\nu^{-1}(p)=\{p',p''\}$.  
Consider the exact sequence
$$0 \ra \cT_{C^{\nu}} \ra \cT_X|_{C^{\nu}} \ra \cN_{C^{\nu}/X} \ra 0.$$
There is an induced extension (see \cite{GHS}\S 2, \cite{AK}\S 25-27):
\begin{equation}
\label{EX1}
0 \ra \cN_{C^{\nu}/X} \ra \nu^*\cN_{C/X} \ra Q \ra 0,
\end{equation}
where $Q$ is torsion, supported at $\{p',p''\}$ with length one
at each point.  We interpret $\nu^*\cN_{C/X}$ in (\ref{EX1})
as the sheaf of sections of $\cN_{C^{\nu}/X}$ with simple poles at $p'$
(resp. $p''$) in the direction $\cT_{C^{\nu}}|_{p''}$ (resp. $\cT_{C^{\nu}}|_{p'}$.)
The elements listed in (\ref{GEN}) may be regarded as free generators for
both $\nu^*\cN_{C^{\nu}/X}|_{p'}$ and $\nu^*\cN_{C^{\nu}/X}|_{p''}$;
the resulting identification expresses the compatibility condition satisfied 
by sections of $\nu^*\cN_{C/X}$ which descend to $\cN_{C/X}$.

\section{K3 surfaces}
\label{sect:k3}

In this section we work over $\bC$.

Let $S_0$ be a smooth proper complex variety and $D\subset S_0$ 
a smooth divisor. There is a natural exact sequence
$$
0\ra \cT_{S_0}\langle -D\rangle \ra \cT_{S_0}\ra \cN_{D/S_0}\ra 0,
$$
where the first term is derivations with logarithmic zeroes along $D$. 
This induces a long exact sequence
$$
\ra H^1(S_0,\cT_{S_0}\langle -D\rangle)\ra  H^1(S_0,\cT_{S_0})\stackrel{\phi}{\ra} H^1(D,\cN_{D/S_0}) \ra
$$
where the second term parametrizes the first order deformations of $S_0$
and the first term parametrizes first order deformations of the pair $(S_0,D)$ (see \cite{kawamata}). 

Let $\pi: \cS\ra \Delta$ be a deformation of $S_0$
over $\Delta=\Spec(\bC[[z]])$, with $S_0=\pi^{-1}(0)$.  
Consider the extension 
\begin{equation} \label{ext1}
0\ra \cN_{D/S_0}\ra \cN_{D/\cS}\ra \cN_{S_0/\cS}|_D\ra 0
\end{equation}
and the corresponding extension class 
$$
\eta\in \Ext^1_{D}(\cN_{S_0/\cS}|_{D}, \cN_{D/S_0})\simeq  \Ext^1_{D}(\cO_{D}, \cN_{D/S_0})=H^1(D,\cN_{D/S_0}).
$$
Let 
$$
\kappa\in H^1(S_0,\cT_{S_0})
$$
be the Kodaira-Spencer class of $d/dz$;  using the isomorphisms
$$H^1(S_0,\cT_{S_0})\simeq \Ext_{S_0}^1(\cO_{S_0},\cT_{S_0}) \quad
\cN_{S_0/\cS}\simeq \cO_{S_0}$$
we see that $\kappa$ classifies the extension 
$$0 \ra \cT_{S_0} \ra \cT_{\cS}|_{S_0} \ra \cN_{S_0/\cS} \ra 0.$$
Restricting this to $D$ and taking quotients by $\cT_D$, we 
obtain extension (\ref{ext1}).  Thus we obtain:

\begin{lemm}
\label{lemm:gen}
The classes $\phi(\kappa)$ and $\eta$ are proportional by a non-zero constant. 
\end{lemm}

\begin{prop}\label{prop:minusone}
Let $\pi: \cS\ra \Delta$ be a deformation of a K3 surface $S_0$,
containing a smooth rational curve $\ell$. 
Assume that the Kodaira-Spencer class $\kappa$ is not 
in the tangent space to the locus of K3 surfaces containing 
some deformation of $\ell$. Then 
$$
\cN_{\ell/\cS} \simeq \cO_{\bP^1}(-1)\oplus \cO_{\bP^1}(-1).
$$
\end{prop}
\begin{proof}
The moduli space of K3 surfaces is smooth, with tangent space $H^1(S_0,\cT_{S_0})$
at $S_0$;  the subspace parametrizing surfaces containing a deformation of $\ell$
is also smooth with tangent space equal to the image of
$$H^1(\cT_{S_0}\langle -\ell \rangle) \ra H^1(\cT_{S_0}).$$
Moreover, $\kappa$ is not in this subspace precisely when
$\phi(\kappa)\neq 0$.  

We have the extension
$$
0\ra \cN_{\ell/S_0}\ra \cN_{\ell/\cS}\ra \cN_{S_0/\cS}|_{\ell}\ra 0
$$
with $\cN_{\ell/S_0}\simeq \cO_{\bP^1}(-2)$ and $\cN_{S_0/\cS}|_{\ell}\simeq \cO_{\bP^1}$.
The extension class $\eta$ is in the group 
$$\Ext_{\ell}^1(\cN_{S_0/\cS}|_{\ell},\cN_{\ell/S_0})\simeq 
H^1(\bP^1,\cO_{\bP^1}(-2))\simeq \bC$$
hence
$$
\cN_{\ell/\cS}= 
               \begin{cases} 
  \cO_{\bP^1}(-1) \oplus \cO_{\bP^1}(-1) & \text{ if } \eta \neq 0\\
\cO_{\bP^1} \oplus \cO_{\bP^1}(-2) & \text{ if } \eta = 0 
               \end{cases}.
$$
Lemma~\ref{lemm:gen} and the hypothesis guarantee
that $\eta \neq 0$.  
\end{proof}

\

For each $g\ge 2$,
let $\cK_{2g-2}$ denote the moduli space of
K3 surfaces $S$ with a primitive polarization $f$ of degree 
$2g-2$.  For each $d>0$,
let $\cD_{g,d}$ denote the moduli space of triples $(S,f,\ell)$, 
where $(S,f)\in \cK_{2g-2}$ and
$$\Lambda_{g,d}:=\left<f,\ell\right>=\bZ f+\bZ \ell \subset \NS(S)$$
is a saturated
sublattice with intersection form
$$
\begin{array}{c|cc}
 & f & \ell \\
\hline
f    & 2g-2 & d \\
\ell & d    & -2
\end{array}.
$$
Let $\iota:\cD_{g,d}\ra \cK_{2g-2}$ be the induced projection, which is finite 
onto its image.  

Each $\cD_{g,d}$ is an open subspace of the moduli spaces 
of lattice-polarized K3 surfaces of 
Nikulin \cite{Nik1} (see \cite{Dol} for detailed discussion).  
Now $\Lambda_{g,d}$ admits a primitive imbedding into the 
cohomology lattice of a K3
surface, which is unique up to conjugation by automorphisms
of the cohomology lattice \cite{Nik2} 1.14.4.  Surjectivity 
of the period map and the global Torelli theorem \cite{K3book}
Expos\'es X,IX guarantee $\cD_{g,d}$
is nonempty and irreducible.  

\begin{prop} \label{prop:existcurve}
Given $(S,f,\ell)\in \cD_{g,d}$ such that $\NS(S)=\left<f,\ell\right>$,
which is the case for a very general such triple.  Then 
there is a smooth
rational curve in $S$ with divisor class $\ell$.  Moreover,
$f$ is very ample if $g>2$ and induces a branched double cover
$S\ra \bP^2$ if $g=2$.
\end{prop}
\begin{proof}
The assertions of the last sentence are applications of results
on linear series of K3 surfaces by Mayer and Saint-Donat:
\begin{enumerate}
\item{complete linear series have no isolated basepoints
\cite{SD} 3.1;}
\item{an ample divisor $f$ has fixed components only when there is an 
smooth curve $E\subset S$ of genus one
and a smooth rational curve $\Gamma\subset S$ so that $f=gE+\Gamma$
\cite{SD} 8.1;}
\item{$f$ is hyperelliptic only when $f^2=2$ (in which
case $S$ is a double cover of $\bP^2$ branched over a sextic)
or there exists
a smooth curve $E\subset S$ of genus one
with $f\cdot E=2$ \cite{SD} 5.2;}
\item{except in the special cases listed above, $f$
is very ample \cite{SD} 6.1.}
\end{enumerate}
Our hypothesis implies that $\NS(S)$ has discriminant
$$\mathrm{Disc}\left<f,\ell\right>=-(f\cdot f)(\ell \cdot \ell)+
(f\cdot \ell)^2=4(g-1)+d^2>4.$$
In case (2), we have
$$\mathrm{Disc} \left<E,\Gamma\right>=-(E\cdot E)(\Gamma \cdot \Gamma)+
(E\cdot \Gamma)^2=1;$$
in the latter case of (3), the discriminant is $4$.  

Recall the description of the monoid
of effective classes of a polarized K3 surface $(S,f)$
\cite{LP} 1.6:  It is generated by
\begin{enumerate}
\item{$D\in \NS(S)$ with $D^2\ge 0$ and $D\cdot f >0$;}
\item{$\Gamma \in \NS(S)$ with $\Gamma^2=-2$ and $\Gamma \cdot f>0$.}
\end{enumerate}
Moreover, indecomposable effective classes $\Gamma$ with $\Gamma^2=-2$
represent smooth rational curves.  Thus $\ell$ is effective, and
represents a smooth rational curve if it is indecomposable.

Assume the contrary, so 
there exists an indecomposable effective class $m$ with
$\ell-m$ effective and nonzero.  
As $f$ is ample, we have $f\cdot \ell>f\cdot m$.  
Since $m$ is indecomposable, $m\cdot m \ge -2$.  
The discriminants satisfy
$$\mathrm{Disc}\left<f,m\right>=-(f\cdot f)(m\cdot m)+(f\cdot m)^2<
-(f\cdot f)(\ell\cdot \ell)+(f\cdot \ell)^2=\mathrm{Disc}\left<f,\ell\right>.$$
However, since $\left<f,m\right> \subset \left<f,\ell\right>$ we have
$$\mathrm{Disc}\left<f,m\right>\ge
\mathrm{Disc}\left<f,\ell\right>,$$
a contradiction.
\end{proof}

\section{Pencils of K3 surfaces} 
\label{sect:k3pencil}
\begin{defi} \label{defn:pencilK3}
A {\em pencil of K3 surfaces} consists of 
\begin{enumerate} 
\item{a projective irreducible threefold $Y \subset \bP^N$;}
\item{a codimension-two linear subspace $\Lambda \subset \bP^N$
meeting $Y$ transversally along a smooth curve $Z$,
so that a generic hyperplane $\Lambda \subset H \subset \bP^N$
intersects $Y$ in a K3 surface.}
\end{enumerate}
Projection from $\Lambda$
$$p_{\Lambda}:\mathrm{Bl}_{\Lambda} \bP^N \ra \bP^1$$
induces a flat morphism
$$\pi:\cS:=\mathrm{Bl}_Z Y \ra \bP^1$$
with generic fiber a K3 surface. 
\end{defi}

Let $E$ denote the exceptional divisor of the blow-up
$\beta:\mathrm{Bl}_Z Y \ra Y$.   
\begin{enumerate}
\item{Since $Z$ is a complete intersection in $Y$ we have
$$E\simeq \bP(\cN_{Z/Y})\simeq \bP(\cO_Z(1) \oplus \cO_Z(1))\simeq \bP^1 \times Z;$$
projection onto the first factor coincides with $\pi|_E$.  Thus
for each $p\in \bP^1$, the fiber $\cS_p=\pi^{-1}(p)$ intersects $E$ 
along $\{p\}\times Z$.}
\item{For each $z\in Z$, the fiber $E_z:=\beta^{-1}(z)$ meets each
member of the pencil in one point and thus yields a section 
$\mathsf s_z:\bP^1 \ra \cS$
of $\pi$.}
\item{The normal bundle of this section is computed by the exact sequence
$$
\begin{array}{ccccccccc}
0 & \ra & \cN_{E_z/E} &\ra  & \cN_{E_z/\cS}                       &\ra & \cN_{E/\cS}|_{E_z} &\ra & 0 \\
0 & \ra & \cO_{\bP^1} & \ra & \cO_{\bP^1} \oplus \cO_{\bP^1}(-1) & \ra & \cO_{\bP^1}(-1) & \ra & 0. 
\end{array} 
$$}
\item{For each smooth fiber $\cS_p$ and each curve $\ell \subset \cS_p$,
the intersection $Z\cap \ell\neq \emptyset$.   Given $z\in Z\cap \ell$,
$E_z$ meets $\ell$ transversally at $z$, i.e.,
the tangent spaces of $\ell$ and $E_z$ at $z$ are transverse.}
\end{enumerate}
 
\

\begin{exam}
Let $Y\subset \bP^N$ be a threefold with isolated singularities and $\omega_Y=\cO_Y(-1)$;
if $Y$ is smooth then it is a Fano threefold.  We can choose $\Lambda$ so that it meets
$Y$ along a smooth curve, which yields a pencil of K3 surfaces as in Definition~\ref{defn:pencilK3}.
\end{exam}

\begin{prop} \label{prop:makepencil}
Fix $2\leq g \leq 10$.  There exists a projective variety
$X\subset \bP^M$ so that 
the generic K3 surface of degree $2g-2$ can be 
realized as an intersection of $X$ with a linear 
subspace of codimension $\dim(X)-2$, and thus as a member
of a pencil of K3 surfaces.  
\end{prop}
\begin{proof}
The generic K3 surface of degree $2$ is a hypersurface 
of degree $6$ in weighted projective space
$\bP(1,1,1,2)$, which has isolated singularities.  Imbed
$$X:=\bP(1,1,1,2) \hookrightarrow \bP^M$$
using weighted-homogeneous forms of degree six.
Smooth hyperplane sections are degree $2$ K3 surfaces.

For $g=3,4,5,6,7,8,9,10$,
the generic K3 surface of degree $2g-2$ is a complete intersection 
in a generalized flag variety \cite{Mukai}.  We extract the
variety $X$ in each case:
\begin{itemize}
\item{$g=3$ quartic hypersurface in $\bP^3$: take $X$ to be
the $4$-fold Veronese
reimbedding of $\bP^3$;}
\item{$g=4$ complete intersection of a quadric and a cubic
in $\bP^4$: take $X$ to be the $3$-fold Veronese
reimbedding of a smooth quadric hypersurface $Q\subset \bP^4$;}
\item{$g=5$ a generic complete intersection of three quadrics: take
$X$ to be the $2$-fold Veronese reimbedding of $\bP^5$;}
\item{$g=6$ a complete intersection of $\Gr(2,5)$ with
a quadric and a codimension-three linear space: take $X$
to be the $2$-fold Veronese reimbedding of the generic codimension-three 
linear section of $\Gr(2,5)$;}
\item{$g=7$ a codimension-eight linear section 
of the isotropic Grassmannian $\mathrm{IGr}(5,10)\subset \bP^{15}$;}
\item{$g=8$ a codimension-six linear section of the Grassmainnian
$\Gr(2,6)\subset \bP^{14}$;}
\item{$g=9$ a codimension-four linear section of the Lagrangian
Grassmannian $\mathrm{LGr}(3,6)\subset \bP^{13}$;}
\item{$g=10$ a codimension-three linear section of the
flag variety of dimension five, associated with the adjoint representation
of $G_2$, imbedded in $\bP^{13}$.}
\end{itemize}

We obtain pencils of K3 surfaces of degree $2g-2$ as follows:  
Taking generic subspaces
$$\Lambda_{N-2} \subset \Lambda_{N} \subset \bP^M$$
of codimensions $\dim(X)-1$ and $\dim(X)-3$ respectively, there is a pencil
$$\Lambda_{N-2} \subset \Lambda_{N-1}(p) \subset \Lambda_N, \quad p\in \bP^1$$ 
of $(N-1)$-dimensional linear spaces.  Taking $Y=\Lambda_N\cap X$
and $Z=\Lambda_{N-2}\cap X$, we obtain a pencil as in 
Definition~\ref{defn:pencilK3}. 
\end{proof}

\section{Deforming sections in pencils}
\label{sect:DSP}
\begin{prop}\label{prop:oneps}
Let $\pi:\cS \ra \bP^1$ be a pencil of K3 surfaces.
Suppose there exists a point $p\in \bP^1$ so that $\cS_p$ is smooth
and contains a smooth rational curve $\ell$ with normal bundle
$\cN_{\ell/\cS}\simeq \cO_{\bP^1}(-1) \oplus \cO_{\bP^1}(-1)$.  
Then there exists a one-parameter family of sections of $\pi$, the closure
of which contains $\ell$.  
\end{prop}
\begin{proof}
Let $Z$ be the base locus of the pencil.
Choose a point $z\in \ell \cap Z \subset \cS_p$ and consider the nodal
curve 
$$C=\ell \cup_z E_z$$
with normalization $\nu:C^{\nu} \ra C$.  Consider the exact sequence
$$0 \ra \cN_{C^{\nu}/\cS} \ra \nu^*\cN_{C/\cS} \ra Q \ra 0$$
with $Q$ a torsion sheaf, of length one at each of the points 
$\{z',z''\}$ over $z$.  

It suffices to show that $\cN_{C/\cS}$ has no higher cohomology and 
admits a global section not mapping to zero in $Q$ (see Section~\ref{sect:deform}):
Then the Hilbert scheme is smooth at $C$ and contains deformations of $C$
smoothing the node $z$ (cf. \cite{HT}, proof of Proposition 24 
and \cite{GHS} \S 2).
Moreover, these deformations meet the generic fiber of $\pi$ at one
point, and thus are sections of $\pi$ not contained in $Z$.  
Finally, since $\ell$ is a component of a degeneration of these 
sections, it lies in the closure of the surface traced out
by these sections.  

We evaluate $\cN_{C/\cS}$ on each component of $C$:  On $E_z$ we
have 
\begin{equation} \label{EX2}
 0 \ra \cN_{E_z/\cS} \ra \cN_{C/\cS}|_{E_z} \ra Q(E_z) \ra 0
\end{equation}
with $Q(E_z)$ of length one at $z$;  furthermore,
$\cN_{E_z/\cS}\simeq \cO_{\bP^1}(-1) \oplus \cO_{\bP^1}$
with the non-negative summand corresponding to the normal direction
to $E_z$ in $E$.  Since the tangent vector to $\ell$ is not contained
in $E$, the $\cO_{\bP^1}$ component is saturated in the extension
(\ref{EX2}), and the $\cO_{\bP^1}(-1)$ component is {\em not} saturated
(\cite{HT} Proposition 23, cf. Sublemma 27).  In particular, we conclude
$$\cN_{C/\cS}|_{E_z} \simeq \cO_{\bP^1} \oplus \cO_{\bP^1}.$$
As for $\ell$, we have 
$\cN_{\ell/\cS}\simeq \cO_{\bP^1}(-1)^{\oplus 2}$ and an extension
$$0 \ra \cN_{\ell/\cS} \ra \cN_{C/\cS}|_{\ell} \ra Q(\ell) \ra 0$$       
with $Q(\ell)$ of length one.  In this case, the only possibility is
$$\cN_{C/\cS}|_{\ell} \simeq \cO_{\bP^1} \oplus \cO_{\bP^1}(-1).$$

We observe that
\begin{enumerate}
\item{$\cN_{C/\cS}|_{E_z}$ is globally generated and has no higher
cohomology;}
\item{$\cN_{C/\cS}|_{\ell}$ has no higher cohomology and each nonzero
global section is nonzero in $Q(\ell)$.}
\end{enumerate}
Since $C$ is obtained by gluing $E_z$ and $\ell$ at a single point, $\cN_{C/\cS}$
also has no higher cohomology \cite{HT}, Lemma 21.  Fixing a nonzero section
$t$ of $\cN_{C/\cS}|_{\ell}$, we can find a section of $\cN_{C/\cS}|_{E_z}$
agreeing with $t$ at $z$.    
\end{proof}

\section{Construction of the examples}
\label{sect:fe}

\begin{theo}
\label{theo:dense2}
Let $\pi:\cS\ra \bP^1$ be a pencil of K3 surfaces
of degree $2g-2$ with base locus
$Z$. Assume that there exists an infinite sequence of positive integers
$$
d_1  < d_2 < \ldots
$$
such that $\bP^1$ intersects 
$\cD_{g,d_j}$ transversally at some point $p_j$ and the
class $\ell_j \in \Lambda_{g,d_j}$ is represented by a smooth
rational curve $\ell_j \subset \cS_{p_j}$.   
Then $\cS$ has a Zariski dense set of sections of $\pi$. 
\end{theo}
\begin{proof}
We have
$$\cN_{\ell_j/\cS}\simeq \cO_{\bP^1}(-1) \oplus \cO_{\bP^1}(-1)$$
by Proposition~\ref{prop:minusone}.  Proposition~\ref{prop:oneps}
gives a one-parameter family of sections
$$\begin{array}{ccc}
\Sigma_j                         & \stackrel{\sigma_j}{\ra} & \cS \\
{\scriptstyle \psi_j} \downarrow \quad &   &     \\
R_j
\end{array},
$$
with $\sigma_j:\psi_j^{-1}(r)\ra \cS$ a section for $r\in R_j$ generic,
and $\sigma_j(\psi_j^{-1}(r_0))\supset \ell_j$ for the distinguished
point $r_0\in R_j$ around which the deformation was produced.  
It follows that the irreducible surface 
$\widetilde{\Sigma}_j:=\sigma_j(\Sigma_j)$ contains $\ell_j$.  

We claim that $\cup_j \widetilde{\Sigma}_j$ is dense in $\cS$.  Let
$\Xi$ denote the closure of all these surfaces.  If $\Xi \subsetneq \cS$
then $\Xi$ is a union of irreducible surfaces, each dominating $\bP^1$.  
Thus the degree $\Xi\cap \cS_p \subset \cS_p$ is bounded;
In particular, $\Xi$ cannot contain $\ell_j$ when $d_j\gg 0$, 
a contradiction.
\end{proof}

\begin{coro}\label{coro:main}
Let $\cS \ra \bP^1$ be a very general pencil of degree $2g-2$,
for $2\le g \le 10$.  
Then the hypothesis of 
Theorem~\ref{theo:dense2} are satisfied, and sections are Zariski dense.
\end{coro}
\begin{proof}
By Proposition~\ref{prop:makepencil}, the moduli space 
$\cK_{2g-2}$ is dominated by 
an open subset of the Grassmannian $\Gr(N,M+1)$ parametrizing the $(N-1)$-dimensional
linear subspsaces in $\bP^M$.  
For $d\gg 0$, $\cD_{g,d}$ 
determines a nonempty divisor in $\Gr(N,M+1)$, which is necessarily ample.
Each pencil determines a $\bP^1 \subset \Gr(N,M+1)$, so the pencils
are parametrized by a suitable Hilbert scheme.  By 
Proposition~\ref{prop:existcurve}, for a generic pencil
$\bP^1$ is transversal to $\cD_{g,d}$ at some
point and the corresponding fiber contains a smooth rational
curve of degree $d$.   Take the countable intersection of these open subsets 
in the space of all pencils;  applying the Baire 
category theorem, we are left with a nonempty dense subset of the space of all
pencils.    
\end{proof}
\begin{rema}\label{rema:uncountable}
Our argument does not preclude the (unlikely) possibility that every pencil defined over
$\overline{\bQ}$ is tangent to $\cD_{g,d}$ for each $d \gg 0$.  
\end{rema}

We complete the argument for Theorem~\ref{theo:one}:
\begin{proof}
Corollary~\ref{coro:main}
gives pencils of K3 surfaces with dense rational points.
Proposition~\ref{prop:makepencil} implies that the generic K3 surface
of degree $2g-2,2\le g \le 10$ arises in such a pencil;  in particular,
K3 surfaces with N\'eron-Severi rank one occur.
The N\'eron-Severi group of the geometric generic fiber
injects in the N\'eron-Severi group of each smooth member
of the pencil;  indeed, the map on homology is an isomorphism.
Thus the N\'eron-Severi group of the geometric generic
fiber is also of rank one. 

For simplicity, we will only discuss non-isotriviality in the
degree four case:  Any pair of quartic surfaces 
$$S_i=\{f_i(w,x,y,z)=0 \}\subset \bP^3, \quad i=1,2$$
are contained in a pencil, e.g., 
$$\cS=\{s_1f_1+s_2f_2=0\} \subset \bP^3 \times \bP^1.$$
This is non-isotrivial whenever $S_1 \not \simeq S_2$, and
completes the proof when $B=\bP^1$.

For an arbitrary smooth complex curve $B$, we can express
$\bC(B)$ as a finite extension of $\bC(\bP^1)$.  Any
K3 surface with dense $\bC(\bP^1)$-rational points 
{\em a fortiori} has dense $\bC(B)$-rational points.
\end{proof}

\section{Limitations}
\label{sect:limit}
In this section, we discuss constraints to deforming sections in
K3 fibrations.  We hope to explain why pencils are a natural source
of examples of K3 surfaces
with potentially-dense rational points.

Let $\pi:\cS \ra B$ be a fibration satisfying the following
\begin{enumerate}
\item{$\cS$ is smooth and the generic fiber is a K3 surface;}
\item{the singular fibers have at worst rational double points;}
\item{$\pi$ is non-isotrivial and projective.}
\end{enumerate}
Let $\mu:B \ra {\overline \cK_{2g-2}}$
be the classifying map into the Baily-Borel compactification of the
corresponding moduli space of polarized K3 surfaces.  Recall that the 
pull-back of the natural polarization on the Baily-Borel compactification is
$\cL:=\pi_*\omega_{\pi}$.  
Since the fibers of $\pi$ have trivial dualizing sheaf, 
$\omega_{\pi}$ is the pull-back of a line bundle from $B$, which 
is necessarily $\cL$ (by the projection formula).  

Suppose we have
a section $\sigma:B\ra \cS$.  Consider the exact sequence
$$0 \ra \cT_{\sigma(B)} \ra \cT_{\cS}|_{\sigma(B)} \ra \cN_{\sigma(B)/\cS} \ra 0.$$
It follows that 
$$c_1(\cN_{\sigma(B)/\cS})=c_1(\cT_{\cS})-c_1(\cT_{\sigma(B)})=-c_1(\omega_{\pi}),$$
whence
$$\deg(\cN_{\sigma(B)/\cS})=-\deg(\mu(B))<0.$$
(For pencils of 
quartic surfaces
the degree is $-1$, the smallest possible value.)
Hence $\cN_{\sigma(B)/\cS}$ cannot be globally generated at the generic point:  If it were
then $\det(\cN_{\sigma(B)/\cS})$ would have a nonvanishing section,
contradicting our degree computation.  
In particular, we cannot deform a section in 
a two-parameter family sweeping out a dense subset of $\cS$.

\section{Projective geometry of pencils of quartic surfaces}
\label{sect:proj}
Corollary~\ref{coro:main} says that a very general pencil of quartic
surfaces has a Zariski-dense set of sections.  Here we interpret
these in terms of classical projective geometry.  

Recall the set-up of Theorem~\ref{theo:dense2}:
$$
\begin{array}{ccc} 
\cS & \stackrel{\pi}{\ra} & \bP^1 \\
{\scriptstyle \beta} \downarrow \quad & & \\
\bP^3 & & 
\end{array}
$$
where $\beta$ is the blowup along the base locus $Z$ of the pencil. 
We analyze the proper transforms in $\bP^3$ 
of the sections of $\pi$
produced by
Proposition~\ref{prop:oneps}.

Let $\ell_1$ be a line in some fiber of $\pi$; $\beta(\ell_1)$ 
is a four-secant line to $Z$. If $z\in Z\cap \beta(\ell_1)$ and 
$E_z:=\beta^{-1}(z)$ then $\beta(\ell_1)\cup_z E_z$ 
deforms to a 3-secant line to $Z$. 
These 3-secants move in a 1-parameter family. 
Similarly, if $\ell_d$ is a smooth rational curve of degree $d$ in some 
fiber of $\pi$ then $\beta(\ell_d)$ is a rational curve of degree $d$
meeting $Z$ in $4d$ points. The nodal curve $\beta(\ell_d)\cup_z E_z$ 
deforms to a 1-parameter family of rational curves in $\bP^3$ of degree $d$ meeting $Z$ 
in $4d-1$ points. 

Essentially, the proof of Proposition~\ref{prop:oneps}
shows that the space of rational curves of degree $d$ in $\bP^3$
meeting $Z$ in at least $4d-1$ points (counted with multiplicities)
has at least one irreducible component of the expected dimension
and the generic such curve meets $Z$ in $4d-1$ points.

\section{Higher-dimensional Calabi-Yau hypersurfaces}
\label{sect:CY}
We have seen in Corollary~\ref{coro:main} that a very general
pencil of quartic surfaces has a Zariski-dense set of sections.
The following generalization is due to Coskun
and Starr:
\begin{theo} \label{theo:dense3}
Let $\pi:\cX\ra \bP^1$ be a very general pencil of hypersurfaces
of degree $n+1$ in $\bP^n$, with $n\ge 4$.  Then sections of
$\pi$ are Zariski dense.
\end{theo}
The first step is the following extension of Proposition~\ref{prop:oneps}
to higher-dimensional pencils (which has an analogous proof):
\begin{prop}\label{prop:onepsnew}
Let $\pi:\cX \ra \bP^1$ be a pencil of hypersurfaces of degree $n+1$
in $\bP^n$.
Suppose there exists a point $p\in \bP^1$ so that $\cX_p$ is smooth
and contains a smooth rational curve $\ell$ with normal bundle
$$\cN_{\ell/\cX}\simeq 
\cO_{\bP^1}^{n-3} \oplus \cO_{\bP^1}(-1) \oplus \cO_{\bP^1}(-1).$$
Then there exists an $(n-2)$-parameter family of sections of $\pi$;
its closure contains the deformations of $\ell$ in $\cX_p$.
\end{prop}

Just as in Theorem~\ref{theo:dense2}, Proposition~\ref{prop:onepsnew}
reduces Theorem~\ref{theo:dense3} to 
the following result, which is a natural extension of
Theorem 1.27 of \cite{Cl} for quintic threefolds:
\begin{theo}
Let $X$ be a very general hypersurface of degree $n+1$ in $\bP^n$,
with $n\ge 4$.
Then there exists an infinite sequence of integers 
$$d_1<d_2<\ldots $$
such that $X$ contains a rational curve $\ell$ of degree $d_j$
with 
$$\cN_{\ell/X}\simeq 
\cO_{\bP^1}^{n-3} \oplus \cO_{\bP^1}(-1) \oplus \cO_{\bP^1}(-1).$$
Furthermore, deformations of these curves are Zariski dense in $X$.
\end{theo}
It suffices
to produce one hypersurface $X$ with the desired properties;
even a singular $X$ will work, as long as it is smooth along each $\ell$.  

The proof proceeds by induction, following \cite{Cl}, which
establishes the base case $n=4$.  Consider a hypersurface 
$X\subset \bP^n$ of degree $n+1$ 
containing the codimension-two linear subspace
$\{x_0=x_n=0\}\subset \bP^n$.  The defining equation for $X$ 
can be written uniquely
$$x_0 f_n(x_0,\ldots,x_{n-1})+x_n g_n(x_0,\ldots,x_n)=0 $$
where $f_n$ and $g_n$ have degree $n$.  
Let $H$ denote the hyperplane $\{x_0=0\}$ and 
$$W=\{f_n(x_0,\ldots,x_{n-1})=0\}\subset  \{x_n=0\}\simeq \bP^{n-1},$$
a hypersurface of degree $n$.  
The singular locus $\mathrm{Sing}(X)$ is defined by the equations
$$f_n+x_0\frac{\partial f_n}{\partial x_0} + 
	x_n\frac{\partial g_n}{\partial x_0}=
x_0\frac{\partial f_n}{\partial x_i}+x_n\frac{\partial g_n}{\partial x_i} = 
x_n\frac{\partial g_n}{\partial x_n}+g_n=0,$$
where $i=1,\ldots,n-1.$  This contains the locus
$$\Sigma:=\{x_0=x_n=f_n=g_n=0\};$$
if this is a smooth complete intersection 
then $\mathrm{Sing}(X)=\Sigma$.  

We make the following generality assumptions on $f_n$:
\begin{itemize}
\item{$W$ satisfies the inductive hypothesis, so
there is an infinite
sequence of integers $\{d_j\}$
so that $W$ contains a rational curve $\ell_j$ of 
degree $d_j$ with normal bundle
$$\cN_{\ell_j/W}\simeq 
\cO_{\bP^1}^{n-4} \oplus \cO_{\bP^1}(-1) \oplus \cO_{\bP^1}(-1).$$}
\item{$H$ intersects $W$ and each $\ell_j$ transversely.}
\end{itemize}
The Bertini theorem allows us to achieve the last condition
for a very general choice of $f_n$.  
We also make assumptions on $g_n$:
\begin{itemize}
\item{$\Sigma$ is a smooth complete intersection, so that 
$\mathrm{Sing}(X)=\Sigma$.}
\item{$\Sigma\cap \ell_j=\emptyset$ for each $j$.}
\end{itemize}
In particular, the $\ell_j$ avoid the singularities of $X$.

We can compute the normal bundle $\cN_{\ell_j/X}$ using
the extension
$$0 \ra \cN_{\ell_j/W} \ra \cN_{\ell_j/X} \ra \cN_{W/X}|_{\ell_j} \ra 0.$$
We have
$$\cN_{W\cup H/X}=\cO_{W\cap H}(+1)$$
and since $H$ is Cartier away from $\Sigma$
$$\cN_{W/X}|_{\ell_j}=\cN_{W\cup H/X}|_{\ell_j}(-H)=\cO_{\ell_j}.$$
The extension
$$0\ra \cO_{\bP^1}^{n-4} \oplus \cO_{\bP^1}(-1)^2 \ra \cN_{\ell_j/X} \ra 
\cO_{\bP^1} \ra 0$$
is necessarily trivial and 
$$\cN_{\ell_j/X}\simeq \cO_{\bP^1}^{n-3} \oplus \cO_{\bP^1}(-1)^2,$$
as desired.  

The deformations of
each $\ell_j$ sweep out a divisor in $W$; the union
of these divisors is dense in $W$.  
It remains to check that the deformations of the $\ell_j$ are dense in 
$X$.  Each $\ell_j$ has deformations in $X$ which are not
contained in $W$.  Indeed, we have shown that $\ell_j$ 
admits infinitesimal deformations 
that are nonzero in $\cN_{W/X}|_{\ell_j}$.  In particular, deformations of
$\ell_j$ sweep out an irreducible divisor $D_j \subset X$ not equal to $W$.  
If the union $\cup_{j=1}^{\infty}D_j$ were a proper subvariety
$\Xi\subsetneq X$, then $\Xi\cap W$ would not be dense in $W$, contradicting the
inductive hypothesis.

\bibliographystyle{plain}
\bibliography{pd}

\end{document}